\documentclass{article}
\usepackage[fleqn]{amsmath}
\usepackage{amssymb,latexsym}
\usepackage{xypic}
\usepackage{url}

\newcommand{\pfaff}{\mathrm{Pfaff}}

\newenvironment{proof}%
{ \medskip
  \rmfamily \noindent
  {\bf Proof:\/} } {\hfill $\Box$ \medskip}

\newenvironment{acknowledgements}
{ \medskip
  \rmfamily \noindent
  {\bf Acknowledgements.\/} } { \bigskip}

\newcommand {\iso}         {\cong}

\let\D=\Delta
\let\s=\sigma
\let\l=\lambda
\newcommand{\sib}{\bar{\sigma}}
\let\ol=\overline

\let\d=\delta
\let\e=\varepsilon
\let\phi=\varphi
\let\sub=\subseteq
\let\bs=\backslash
\let\inc=\hookrightarrow
\let\ck=\check

\newcommand{\xig}{\boldsymbol{\xi}}
\newcommand{\zetag}{\boldsymbol{\zeta}}
\newcommand{\alphag}{\boldsymbol{\alpha}}
\newcommand{\taug}{\boldsymbol{\tau}}

\newcommand{\Dcal}{\mathcal{D}}

\newcommand{\Z}{\mathbb{Z}}
\newcommand{\x}{\mathbf{x}}
\newcommand{\m}{\mathbf{m}}
\newcommand{\y}{\mathbf{y}}

\let\mpty=\varnothing

\newcommand{\R}{\ensuremath{\mathbb{R}}}
\newcommand{\N}{\ensuremath{\mathbb{N}}}
\newcommand{\C}{\ensuremath{\mathbb{C}}}

\renewcommand{\S} {\ensuremath{\mathcal{S}}} 
\renewcommand{\P} {\ensuremath{\mathcal{P}}}
\newcommand{\Q} {\ensuremath{\mathcal{Q}}}
\newcommand{\U} {\ensuremath{\mathcal{U}}}

\newcommand{\Frac}[2]{\frac{\displaystyle{#1}}{\displaystyle{#2}}}

\newcommand{\junk}[1]{}

\newtheorem{theorem}{Theorem}
\newtheorem{lemma}[theorem]{Lemma}
\newtheorem{corollary}[theorem]{Corollary}
\newtheorem{definition}[theorem]{Definition}
\newtheorem{proposition}[theorem]{Proposition}
\newtheorem{example}[theorem]{Example}
\newtheorem{rem}[theorem]{Remark}

\begin{document}

\title{Topology of definable Hausdorff limits}
\author{Thierry Zell}

\maketitle

\begin{abstract}
Let $A\sub \R^{n+r}$ be a set definable in an o-minimal expansion $\S$
of the real field, $A' \sub \R^r$ be its projection, and assume that
the non-empty fibers $A_a \sub \R^n$ are compact for all $a \in A'$
and uniformly bounded, {\em i.e.} all fibers are contained in a ball
of fixed radius $B(0,R).$
If $L$ is the Hausdorff limit of a sequence of fibers $A_{a_i},$ we
give an upper-bound for the Betti numbers $b_k(L)$ in terms of
definable sets explicitly constructed from a fiber $A_a.$ In
particular, this allows to establish effective complexity bounds in the
semialgebraic case and in the Pfaffian case. In the Pfaffian setting,
Gabrielov introduced the {\em relative closure} to construct the
o-minimal structure $\S_\pfaff$ generated by Pfaffian functions in a
way that is adapted to complexity problems.
Our results can be used to estimate the Betti numbers of a relative closure
$(X,Y)_0$ in the special case where $Y=\mpty.$
\end{abstract}

\section*{Introduction}

Let us consider a bounded subset $A\sub \R^{n+r}$ which is definable
in an o-minimal expansion $\S$ of the real field (the reader can refer
to~\cite{coste} or~\cite{vdd98} for definitions).  Let $A'$ be the
canonical projection of $A$ in $\R^r,$ and for all $a \in A',$ we
define the fiber $A_a$ as $A_a=\{\x \in \R^n \mid (\x,a)\in A\}.$ Assume
that these fibers are compact for all $a\in A'.$ Note that since we
assumed that $A$ was bounded, the fibers $A_a$ are all contained in a
ball $B(0,R)$ for some $R>0.$ 
Recall that for compact subsets $A$ and $B$ of $\R^n,$ we can define
the {\em Hausdorff distance} between $A$ and $B$ as 
\begin{equation*}
d_H(A,B)=\max_{\x\in A}\min_{\y\in B} |\x-\y|
+\max_{\y\in B} \min_{\x \in A} |\x-\y|.
\end{equation*}
The Hausdorff distance gives the space $\mathcal{K}_n$ of compact
subsets of $\R^n$ a metric space structure.

\medskip

If $(a_i)$ is a sequence in $A',$ and $L$ is a compact subset of
$\R^n$ such that the limit of the sequence $d_H(A_{a_i},L)$ is zero,
we call $L$ the {\em Hausdorff limit} of the sequence $A_{a_i}.$ It is
a well-established fact that when $A$ is definable in an o-minimal
structure $\S,$ then the Hausdorff limit $L$ is also definable in
$\S:$ it was first proved by Br\"ocker~\cite{brocker} in the algebraic
case; in the general case, it follows from the definability of types
that was first proved by Marker and Steinhorn~\cite{marker}, and later by
Pillay~\cite{pillay}.  Recently, direct proofs were suggested, one
using model theoretic arguments by van~den~Dries~\cite{vdd:H} and a
purely geometric one by Lion and Speissegger~\cite{ls:Hausdorff}.

\subsection*{Main result}

In this paper, we investigate how the topology of the Hausdorff limit
can be related to the topology of the fibers $A_a$ and their Cartesian
powers. To do so, we need to introduce for any integer $p$ a distance
function $\rho_p$ on $(p+1)$-tuples $(\x_0, \ldots, \x_p)$ of points in
$\R^n$ by
\begin{equation}\label{eq:rho}
\rho_p(\x_0, \ldots,\x_p)=\sum_{0 \leq i < j \leq p} |\x_i-\x_j|^2;
\end{equation}
(where $|\x|$ is the Euclidean distance in $\R^n$).  
The {\em expanded $p$-th diagonal} of $A_a$ is defined for all $\d>0$ by
\begin{equation}\label{eq:diag}
D^p_a(\d)=\{(\x_0, \ldots, \x_p) \in (A_a)^{p+1} \mid \rho_p(\x_0,
\ldots, \x_p)\leq \d\}.
\end{equation}

\bigskip

Let $b_k(L)$ denote the $k$-th Betti number of $L,$ by which we mean
the rank of the singular homology group $H_k(L,\Z)$.  Our main result
is the following upper-bound.

\begin{theorem}
\label{th:main}
Let $A \sub \R^{n+r}$ be a bounded definable set with compact fibers
and $L$ be the Hausdorff limit of some sequence $A_{a_i}.$ Then, there
exists $a\in A'$ and $\d>0$ such that for any integer $k,$ we have
\begin{equation}\label{eq:bkL}
b_k(L) \leq \sum_{p+q=k} b_q(D^p_a(\d));
\end{equation}
where the set $D^p_a(\d)$ is the expanded $p$-th diagonal defined
in~\eqref{eq:diag}. 
\end{theorem}

The proof of this theorem relies on the construction of a continuous
surjection from some fiber $A_a$ to $L,$ and the use of the
spectral sequence associated to such a surjection that was already
used in~\cite{gvz}. The spectral sequence alone does not provide
directly an estimate in terms of the topology of explicit sets such as
the sets $D^p_a(\d):$ the bound~\eqref{eq:bkL} is finally obtained
after an approximation process.

\bigskip

Thus, Theorem~\ref{th:main} allows to estimate the Betti numbers of
$L,$ -- which is a definable set, but not obviously so, -- in terms of
the Betti numbers of the sets $D^p_a(\d)$ which are not only clearly
definable, but also easy to describe from a formula defining $A.$ In
particular, if $A$ is defined by a quantifier-free formula, the sets
$D^p_a(\d)$ can also be described without quantifiers.  In the
semialgebraic setting, this allows us to give good effective bounds on
the Betti numbers of Hausdorff limits, since the terms
$b_q(D^p_a(\d))$ are easy to bound. When $A$ is given by a
quantifier-free formula (see Definition~\ref{def:qf}), we obtain the
following estimates.

\begin{corollary}\label{cor:algebraic}
Let $A \sub \R^{n+r}$ be a bounded semialgebraic set with compact
fibers and defined by a quantifier-free sign condition
$\Phi(\x,a)$ on a family 
\begin{equation*}
\P=\{p_1(\x,a), \ldots, p_s(\x,a)\}
\end{equation*}
of polynomials (where $\x \in \R^n$ and $a \in \R^r$).  Let $L$ be
the Hausdorff limit of some sequence of fibers $A_{a_i}.$ If
$\deg_\x(p_j)\leq d$ for all $1 \leq j \leq s,$ we have for any
integer $k,$
\begin{equation*}
b_k(L)\leq O(k^2s^2d)^{(k+1)n}.
\end{equation*}
\end{corollary}

In particular, the Betti numbers of the Hausdorff limit do not depend
on the degrees in $a$ of the polynomials of $\P.$

\subsection*{Application to the Pfaffian structure} 

The present work was motivated by the case where $\S$ is the o-minimal
structure generated by Pfaffian functions.  This class of
real-analytic functions was introduced by Khovanskii~\cite{k91}; it
contains many of the so-called {\em tame} functions that can appear in
applications, such as real elementary functions or Liouville
functions. They are also the basis for the theory of {\em fewnomials},
the study of the behaviour of real polynomials in terms of the number
of monomials that appear with a non-zero coefficient.  (See
section~\ref{sec:estimates} for definitions.)  Wilkie proved
in~\cite{wilkie99} that Pfaffian functions generate an o-minimal
structure $\S_\pfaff$; this result was generalized in~\cite{KaMcI,lr98,s99}.

\medskip

Pfaffian functions are endowed with a natural notion of complexity, or
{\em format} (see Definition~\ref{df:qformat}), which is a tuple of
integers that control their behaviour. This translates easily into a
notion of format for sets defined by quantifier-free formulas (called
{\em semi-Pfaffian} sets). However, the structure $\S_\pfaff$ contains
sets that cannot be defined by a quantifier-free sign condition on
Pfaffian functions.
In~\cite{g00}, Gabrielov gave an alternative to Wilkie's construction
of $\S_\pfaff,$ showing that definable sets could be constructed by
allowing the operation of {\em relative closure} on 1-parameter
couples of semi-Pfaffian sets.  The object of this construction was to
extend the notion of format to all Pfaffian sets, and use it to
generalize the quantitative results already known for semi- and
sub-Pfaffian sets (see the survey~\cite{gv:survey} and references).

\bigskip

The relative closure is defined as follows: we consider $X$ and $Y$
semi-Pfaffian subsets of $\R^n\times \R_+$ as families of
semi-Pfaffian subsets of $\R^n$ depending on a parameter $\l>0.$ When
the couple $(X,Y)$ verifies additional properties (see~\cite{g00} for
details), the relative closure of $(X,Y)$ is defined as $(X,Y)_0=\{\x
\in \R^n \mid (\x,0) \in \ol{X}\bs\ol{Y}\},$ where $\ol{X}$ is the
topological closure of $X.$ In the special case where $Y=\mpty,$ we
denote the relative closure by $X_0.$ When $Y$ is empty, the
restrictions put on couples imply that the fibers $X_\l$ are compact,
and $X_0$ is then simply the Hausdorff limit of $X_\l$ as $\l$ goes to
zero.

\bigskip

Theorem~\ref{th:main} is applicable in this special case, and
effective estimates can be derived as in the algebraic case, since as
in the case of Corollary~\ref{cor:algebraic}, the set $D^p_a(\d)$ is
given by a quantifier-free formula. We obtain

\begin{corollary}\label{cor:rcbnd}
Let $X\sub \R^n \times \R_+$ be a bounded semi-Pfaffian set such that the
fiber $X_\l$ is compact for all $\l>0,$ and let $X_0$ be the relative
closure of $X.$ If for $\l$ small enough, the format of $X_\l$ is bounded
component-wise by $(n,\ell, \alpha, \beta, s),$ we have for any
integer $k,$
\begin{equation}\label{eq:H-betti}
b_k(X_0) \leq  2^{\ell^2(k+1)^2/2} \, s^{2n(k+1)} \,
O(kn(\alpha+\beta))^{(k+1)(n+\ell)}.
\end{equation}
\end{corollary}

\begin{rem}
Theorem~\ref{th:main} can also be used to derive estimates on the
Betti numbers of a general relative closure $(X,Y)_0$, where $Y$ is not
empty (and thus, $(X,Y)_0$ is not necessarily compact). 
This fact is proved in~\cite{z:these}, and good upper-bounds for that
case be the subject of a separate paper.
\end{rem}

\subsection*{Organization of the paper}

The rest of the paper is organized as follows: in
section~\ref{sec:outline}, we reduce the problem of the Hausdorff
limit of a sequence in a definable family $A \sub \R^{n+r}$ to the
case of the Hausdorff limit $X_0$ of a 1-parameter family $X \sub \R^n
\times \R_+$ when the parameter $\l$ goes to zero.  We then describe
the ingredients of the proof of Theorem~\ref{th:main} for that case: we
need to construct a family of continuous surjections $f^\l: X_\l \to
X_0.$ Using the spectral sequence associated to such a surjection, we
can estimate the Betti numbers of $X_0$ in terms of the Betti numbers
of the fibered products of $X_\l.$ Such fibered products need then to
be approximated to obtain an estimate in terms of the Betti numbers of
expanded diagonals $D^p_\l(\d).$

\medskip

In section~\ref{sec:triang}, the family $f^\l$ is constructed using
definable triangulations of functions. We prove two important
properties of this family: $f^\l$ is close to identity when $\l$ goes
to zero and for any $\l' \neq \l,$ we can obtain $f^\l$ by composing
$f^{\l'}$ with an homeomorphism $h: X_\l \to X_{\l'}$ (see
Proposition~\ref{prop:fl}).

\medskip

Section~\ref{sec:approx} is devoted to the topological approximations
that lead to Theorem~\ref{th:main}. In section~\ref{sec:estimates},
the algebraic and Pfaffian complexity estimates
(Corollary~\ref{cor:algebraic} and Corollary~\ref{cor:rcbnd}) are
proved. The section contains also all the relevant background material
on Pfaffian functions and on Betti numbers of quantifier-free
formulas.

\section{Reduction to one parameter and strategy}
\label{sec:outline}

In this section, we show how, using the results of Lion and
Speissegger on definability of Hausdorff limits~\cite{ls:Hausdorff},
we can reduce the general case of a Hausdorff limit that occurs in a
family with $r$ parameters to the case where $r=1.$ 

\bigskip

Fix $\S$ an o-minimal expansion of the real field
(see~\cite{coste,vdd98}). Let $A\sub \R^{n+r}$ be a bounded
definable set with compact fibers, and $L$ be the Hausdorff limit of a
sequence of fibers of $A.$ We assume of course that $L$ is not already
a fiber of $A,$ since Theorem~\ref{th:main} is trivial in this case.
Since sequences of parameters $(a_i)$ in $A'$ are not definable in
$\S,$ it is difficult to handle Hausdorff limits directly. To avoid
this problem, Lion and Speissegger constructed in~\cite{ls:Hausdorff}
a new family $B$ to model the Hausdorff limits of fibers of $A.$ 
The main result they prove is the following.

\begin{theorem}[\cite{ls:Hausdorff}]\label{th:ls}
If $A\sub \R^{n+r}$ is a bounded definable set with compact fibers,
there exists $R \geq r$ and a definable compact $B \sub \R^{n+R}$ such
that, if $B'$ is the projection of $B$ on $\R^R,$ the following
properties hold.
\begin{itemize}
\item[(H1)] For every $a\in A',$ there is a $b \in B'$ such that $A_a=B_b.$;
\item[(H2)] for every sequence $(b_i)$ in $B'$ such that $\lim b_i=b^*,$
the Hausdorff limit of $B_{b_i}$ exists and equals $B_{b^*};$
\item[(H3)] $\dim B'=\dim A'$ and $\dim \{b \in B' \mid \forall a \in
A', B_b\neq A_a\}< \dim A'.$
\end{itemize}
\end{theorem}

The proof of this result is quite technical, and involves representing
the fibers of $A$ and their possible Hausdorff limits in terms of
integral manifolds of some distributions that depend only on $A.$

\medskip

Using Theorem~\ref{th:ls}, we can obtain $L$ as a limit of a 1-parameter
family by the following proposition.

\begin{proposition}~\label{prop:red}
Let $A\sub \R^{n+r}$ be a bounded definable set with compact fibers,
and $L$ be the Hausdorff limit of a sequence $A_{a_i}.$ Then, there
exists a definable family $X \sub \R^n\times (0,1)$ such that the
following holds.
\begin{itemize}
\item[(X1)] For every $\l\in (0,1),$ there exists $a(\l)\in A'$ such that
  $X_\l=A_{a(\l)}.$ 
\item[(X2)] $L$ is the Hausdorff limit $X_0$ of $X_\l$ when $\l$ goes to zero.
\end{itemize}
\end{proposition}

\begin{proof}
Let $B$ be the set described in Theorem~\ref{th:ls}. By property (H1),
the set of parameters $B'$ contains a sequence $b_i$ such that
$A_{a_i}=B_{b_i}$ for all $i.$ Since $B'$ is compact, we can assume by
taking a subsequence that $b_i$ converges to some $b^* \in B'.$ By
property (H2), we must have $L=B_{b^*},$ since the Hausdorff limit
is unique.
Since $A_{a_i}=B_{b_i}$ for all $i,$ the point $b^*$ is in the closure
of the definable set 
\begin{equation*}
C=\{b \in B' \mid \exists a \in A', A_a=B_b\}.
\end{equation*}
By the curve selection lemma \cite[Chapter~6, Corollary~1.5]{vdd98},
there exists a definable curve $\gamma: (0,1) \to C$ such that
$\lim_{\l \to 0} \gamma(\l)=b^*.$ 
Consider the definable family $X$ given by
\begin{equation*}
X=\{(\x,\l) \in \R^n \times (0,1) \mid \x \in B_{\gamma(\l)}\}.
\end{equation*}
Since $\gamma(\l)\in C$ for all $\l \in (0,1),$ there exists for each
$\l$ a point $a(\l) \in A'$ such that $ B_{\gamma(\l)}=A_{a(\l)},$ so
property (X1) holds. Moreover, property (H2) in Theorem~\ref{th:ls}
guarantees that the Hausdorff limit of $B_{\gamma(\l)}$ when $\l$ goes
to zero is $B_{b^*},$ and since $B_{b^*}=L$ by construction, (X2)
holds too.
\end{proof}

\bigskip

Throughout the rest of this paper, we will assume we are given $X$ as
in Proposition~\ref{prop:red}, with Hausdorff limit $X_0$ when $\l$
goes to zero. The strategy will be the following.
Using the fact that $\ol{X}$ is definable and compact, and thus that
the projection $\pi$ of $\ol{X}$ on the $\l$-axis can be triangulated,
we will construct in section~\ref{sec:triang} a family of continuous
surjections $f^\l: X_\l \to X_0$ defined for small values of $\l.$
Since $X_\l$ and $X_0$ are both compact, the surjection $f^\l$ is {\em
closed}, and we can apply the following theorem~\cite[Theorem~1]{gvz}.

\begin{theorem}\label{th:ss}
Let $f: X \to X_0$ be a closed continuous surjective map definable in
an o-minimal structure.  For all integers $k,$ the following inequality
holds.
\begin{equation}\label{eq:th0}
b_k(X_0) \leq \sum_{p+q=k} b_q(W^p);
\end{equation}
where $W^p$ is the $(p+1)$-fold fibered product of $X;$
\begin{equation}\label{eq:Wp-def}
W^p=\{(\x_0,\ldots, \x_p) \in X^{p+1} \mid f(\x_0)=\cdots=f(\x_p)\}.
\end{equation}
\end{theorem}

Thus, the existence of $f^\l$ for a fixed $\l>0$ gives estimates on
the Betti numbers of $X_0$ in terms of the Betti numbers of the
definable sets $W^p_\l$ (obtained by taking $f$ to be $f^\l$ in
Theorem~\ref{th:ss}). However, there is no explicit description of
these fibered products in the general case,
so this fact alone is not sufficient to establish effective
upper-bounds in the algebraic and Pfaffian case.

\bigskip

Section~\ref{sec:approx} is devoted to refining the estimate given by
the spectral sequence to finally obtain
Theorem~\ref{th:main}, which gives an estimate for the Betti
numbers of $X_0$ in terms of definable sets that are described in a
completely explicit way: the expanded diagonals $D^p_\l(\d).$ The
result is achieved by showing that for suitable values of $\d$ and
$\l,$ the fibered product $W^p_\l$ is included in the expanded
diagonal $D^p_\l(\d),$ and that this inclusion induces an isomorphism
between the corresponding homology groups
(Proposition~\ref{prop:approx}).

\section{Construction of a family of surjections}
\label{sec:triang}

The setting for this section is the following: we consider a definable
family $X \sub \R^n\times (0,1)$ such that the fiber $X_\l$ is compact
for all $\l \in (0,1)$ and such that this family has a Hausdorff limit
$X_0$ when $\l$ goes to zero. As announced in the previous section, we
will construct for small values of $\l$ a family of continuous
surjections $f^\l: X_\l \to X_0$ that are close to identity. More
precisely, we will prove the following result.

\begin{proposition}\label{prop:fl}
Let $X$ be a definable family as above. 
There exists $\l_0>0$ and a family of definable continuous surjections 
$f^\l: X_\l \to X_0,$ defined for all $\l \in (0,\l_0)$ such that 
\begin{equation}\label{eq:norm2}
\lim_{\l \to 0} \, \max_{\x \in X_\l} |\x-f^\l(\x)|=0.
\end{equation}
Moreover, this family $f^\l$ verifies the following property: 
for all $0<\l'<\l<\l_0,$ there exists a (uniformly)
definable homeomorphism $h: X_\l \to X_{\l'}$ such that 
for all $\x\in X_\l,$ we have $f^\l(\x) =f^{\l'}(h(\x)).$ 
\end{proposition}

\subsection{Triangulation of the projection on $\l$}

Since the terminology concerning simplexes is somewhat variable, we
will now state the precise definitions we will be using. These
definitions will follow~\cite{coste} rather than~\cite{vdd98}.

\begin{definition}
If $\alphag_0, \ldots, \alphag_d$ are affine-independent points in
$\R^n,$ the {\em closed simplex} $\sib=[\alphag_0, \ldots, \alphag_d]$
is the subset of $\R^n$ defined by
\begin{equation}\label{eq:closed-simplex}
\sib=\left\{\sum_{i=0}^d w_i\,\alphag_i \mid \sum_{i=0}^d w_i=1, w_1\geq 0,
\ldots, w_d\geq 0\right\}.
\end{equation}
A function $g: \sib \to \R$ is {\em affine} if it satisfies the equality
\begin{equation}\label{eq:affine}
g\left( \sum_{i=0}^d w_i\, \alphag_i \right)= \sum_{i=0}^d w_i \,g(\alphag_i);
\end{equation}
for any $w_0, \ldots ,w_d$ as in~\eqref{eq:closed-simplex}.
\end{definition}

A {\em face} of $\sib$ is any closed simplex obtained from a non-empty
subset of $\alphag_0, \ldots, \alphag_d.$ The {\em open} simplex
$\s=(\alphag_0, \ldots, \alphag_d)$ is the subset of points
$\sum_{i=0}^d w_i\,\alphag_i$ in $\sib$ for which $w_i>0$ for all $0
\leq i \leq d.$

\medskip

\begin{definition}
A {\em (finite\/) simplicial complex} $K$ of $\R^n$ is a finite
collection $\{\sib_1, \ldots, \sib_k\}$ of closed simplexes that is
closed under taking faces, and such that $\sib_i \cap \sib_j$ is a
common face of $\sib_i$ and $\sib_j$ for any $1\leq i,j \leq k.$ The
{\em geometric realization} of $K$ is the subset of $\R^n$ defined by
$|K|=\sib_1 \cup \cdots \cup \sib_k.$
\end{definition}

\bigskip

Throughout section~\ref{sec:triang}, we will denote by $\pi: \ol{X}
\to \R$ the projection on the $\l$-coordinate. Note that $\ol{X}$ is
definable (since it is the closure of a definable set, see for
instance~\cite[Proposition~1.12]{coste}), so $\pi$ is definable too
({\em i.e.} the graph of $\pi$ is a definable set).  Since $\ol{X}$ is
also compact, the triangulation theorem for definable functions
(see~\cite{coste:few, coste}) allows us to assume without loss of
generality that $\ol{X}$ is the geometric realization of a simplicial
complex $K$ and that the map $\pi$ is affine on each simplex $\sib$ of
$K.$ (Note that in general, this requires a linear change of
coordinates in the fibers, but this does not affect our results.)
Moreover, we'll identify $X_0$ with $X_0\times\{0\}\sub X$ and we'll
assume that the triangulation has been refined so that it is
compatible with $X_0,$ {\em i.e.}  $X_0$ is the union of open
simplices of the triangulation.

\bigskip

In the present section, we will need to consider points in both the
total space $\ol{X} \sub \R^n\times [0,1]$ and points in fibers
$X_\l,$ which are by definition subspaces of $\R^n.$ To avoid
ambiguities, we will use the following convention: bold Greek letters
such as $\xig$ will denote points in the total space, whereas bold
Roman letters will be used to denote the points in fibers $X_\l.$
Thus, if $\pi(\xig)=\l>0,$ we have $\xig=(\x,\l),$ where $\x\in X_\l.$

\begin{definition}
The {\em star} $S$ of $X_0$ in $X$ is the union of $X_0$ with all the
open simplices $(\alphag_0, \ldots, \alphag_d)$ that have at least one
vertex $\alphag_i$ in $X_0.$
\end{definition}

\subsection{Constructing a retraction}

Let us define the retraction $F:S \to X_0$ as follows. If $\xig \in
X_0,$ we let $F(\xig)=\xig.$ If $\xig$ belongs to some open simplex
$\s=(\alphag_0, \ldots, \alphag_d),$ where $\alphag_0, \ldots,
\alphag_d$ are vertices such that $\alphag_0, \ldots, \alphag_k$ are
in $X_0$ and $\alphag_{k+1}, \ldots, \alphag_d$ are not in $X_0,$ for
some $k$ with $0\leq k<d,$ define $F$ on $\s$ by
\begin{equation}\label{eq:F}
F\left(\sum_{i=0}^d w_i\,\alphag_i\right)=
\frac{1}{\sum_{i=0}^k w_i} \quad \sum_{i=0}^k w_i\,\alphag_i.
\end{equation}
If $\xig=(\x,\l)\in S$ with $\l>0,$ we'll denote by  $\D(\xig)$  the
intersection between the line through $\xig$ and $F(\xig)$ and the
unique open simplex $\s$ containing $\xig.$
\bigskip

\begin{proposition}\label{prop:F}
The above definition gives a continuous retraction $F: S \to X_0$ that
verifies: for all $\xig=(\x,\l)\in S$ with $\l>0$ and all $\zetag \in
\D(\xig),$  we have $F(\zetag)=F(\xig)$.
\end{proposition}

\begin{proof}
Let $\s=(\alphag_0, \ldots, \alphag_d)$ be an open simplex, with
$\alphag_0, \ldots, \alphag_k$ in $X_0$ and $\alphag_{k+1}, \ldots,
\alphag_d$ not in $X_0,$ for some $k$ with $0\leq k<d,$ so that $\s
\sub S.$ Fix $\xig=\sum_{i=0}^d w_i \alphag_i$ in $\s,$ and let
$s=\sum_{i=0}^k w_i.$ Since all the weights $w_i$ are positive, the
inequality $0\leq k<d$ implies that $0<s<1.$ Thus, the
formula~\eqref{eq:F} clearly defines a continuous function from $\s$
to $X_0.$ Moreover, if $\s'=(\alphag_{i_0}, \ldots, \alphag_{i_e})$ is
a face of $\s$ with at least one $0\leq j\leq e$ such that $i_j\leq k$
(so that $\s' \sub S$), it is clear that the expression~\eqref{eq:F}
extends $F$ continuously to $\s'.$

\medskip

If $\s$ and $\xig$ are as above and $\zetag=t \xig+ (1-t) F(\xig)$ is a
point on $\D(\xig),$ we will now show that $F(\xig)=F(\zetag).$ 
We have $\zetag=\sum_{i=0}^k w'_i \, \alphag_i,$ where
\begin{equation*}
w'_i=
\begin{cases}
t w_i + (1-t)\Frac{w_i}{s}
& \hbox{ if } 0 \leq i \leq k; \\
t w_i & \hbox{ if } k+1 \leq i \leq d.
\end{cases}
\end{equation*}
To prove that $F(\xig)=F(\zetag),$ we must prove that for all $0 \leq i \leq
k,$ 
\begin{equation*}
\frac{w_i}{\sum_{j=0}^k w_j}=\frac{w'_i}{\sum_{j=0}^k w'_j}.
\end{equation*}

Cross-multiplying, we get the following quantities.
\begin{equation}\label{cross1}
w_i\, \sum_{j=0}^k w'_j= 
w_i\, \sum_{j=0}^k \left(t w_j+ (1-t)\Frac{w_j}{s} \right)
=w_i \left(1-t+t s\right);
\end{equation}
and 
\begin{equation}\label{cross2}
w'_i\, \sum_{j=0}^k w_j= \left(t w_i +
(1-t)\Frac{w_i}{s}\right)s
=(t s + (1-t))w_i.
\end{equation}
The two final expressions in~\eqref{cross1} and~\eqref{cross2} are clearly
equal, so $F(\zetag)=F(\xig)$ for any $\zetag \in \D(\xig).$ 
\end{proof}

\begin{definition}\label{df:fl}
Let $\l_0=\min\{\pi(\alphag)\mid \alphag \hbox{ is a vertex of }X,\  \alphag
\not \in X_0\}.$ For any $\l \in (0,\l_0),$ we define $f^\l: X_\l \to
X_0$ by $f^\l(\x)=F(\x,\l).$
\end{definition}

Since $\pi(\alphag)=0$ can only happen if $\alphag \in X_0,$ it
follows that $\l_0>0.$
Note also that if $\xig=(\x,\l)$ is not in $S,$ then we must have
$\l\geq \l_0.$ Indeed, if $\xig \not\in S,$ it belongs to an open
simplex $\s$ of the form $(\alphag_0, \ldots, \alphag_d)$ such that
none of the vertices is in $X_0.$ Thus, $\pi(\alphag_i) \geq \l_0$ for
all $i,$ and since $\pi$ is affine on $\s,$ we must have $\l=\pi(\xig)
\geq \l_0$ too. Thus, for any fixed $\l\in(0,\l_0),$ we have
$\{(\x,\l) \mid \x\in X_\l\} \sub S,$ so $f^\l$ is well-defined for
$\l\in(0,\l_0).$ Since $F$ is continuous, $f^\l$ is continuous too.

\subsection{Properties of the maps $f^\l$}

We must still show that the family of mappings $f^\l$ has all the
properties described in Proposition~\ref{prop:fl}: $f^\l$ is surjective
(Lemma~\ref{lem:surjective}), close to identity
(Proposition~\ref{prop:normF}), and $f^\l$ can be obtained from
$f^{\l'}$ by composing on the right by a homeomorphism $h:X_\l \to
X_{\l'}$ (Proposition~\ref{prop:homeo}).

\begin{lemma}\label{lem:surjective}
For all $\l \in (0,\l_0),$ the map $f^\l$ is surjective.
\end{lemma}

\begin{proof}
Let $\zetag \in X_0.$ Then, there exists a unique set of vertices
$\{\alphag_0, \ldots, \alphag_k\}$ such that $\zetag$ belongs to the
open simplex $(\alphag_0, \ldots, \alphag_k);$ let $v_0, \ldots, v_k$
be the corresponding weights, so that $\zetag=\sum_{i=0}^k v_i
\alphag_i.$ There must be vertices $\alphag_{k+1}, \ldots, \alphag_d$
such that the open simplex $\s=(\alphag_0, \ldots, \alphag_d)$ is in
$X,$ otherwise $\zetag$ could not be approximated by points of $X_\l$
for $\l>0.$

\medskip

Let $\xig=\sum_{i=0}^d w_i \alphag_i$ where $w_i=v_i/2$ for $0 \leq i
\leq k$ and $w_{k+1}, \ldots, w_d$ are arbitrarily chosen positive
numbers so that $\sum_{i=0}^d w_i=1.$ By choice of $w_0, \ldots, w_k,$
we have $\sum_{i=0}^k w_i=1/2,$ and thus $F(\xig)=\zetag.$ Moreover,
if $\D(\xig)$ is as defined in Proposition~\ref{prop:F}, there must be
a point $\taug\in \D(\xig)$ such that $\pi(\taug)=\l.$ Indeed, if we
parameterize the line between $\xig$ and $\zetag$ by $\{(1-t) \zetag +
t \xig \mid t \in \R\},$ the endpoints of $\D(\xig)$ are obtained for
$t=0$ and $t=2,$ which give respectively the points $\zetag$ and
$\zetag'=2\sum_{i=k+1}^d w_i \alphag_i.$ We have $\pi(\zetag)=0,$ and
since $\zetag'$ is not in $S,$ we have $\pi(\zetag')\geq \l_0.$ Since
by restriction $\pi$ is affine on $\D(\xig),$ $\pi(\taug)$ takes all
the values in the interval $(0, \pi(\zetag'))$ when $\taug$ runs
through $\D(\xig).$ In particular, if $\l < \l_0$ there exists $\taug
\in \D(\xig)$ with $\pi(\taug)=\l.$ By Proposition~\ref{prop:F}, we
must have $F(\taug)=F(\xig)$ and since $F(\xig)=\zetag,$ this proves
that $f^\l$ is surjective.
\end{proof}

\begin{proposition}\label{prop:normF}
For $f^\l$ as in Definition~\ref{df:fl}, we have
\begin{equation}\label{eq:limFl}
\lim_{\l\to 0} \, \max_{\x\in X_\l} |\x-f^\l(\x)| =0.
\end{equation}
\end{proposition}

\begin{proof}
Let $\s=(\alphag_0, \ldots, \alphag_d)$ be an open simplex where
$\alphag_0, \ldots, \alphag_k$ are in $X_0$ and $\alphag_{k+1},
\ldots, \alphag_d$ are not in $X_0,$ where $0\leq k<d.$ Fix $\xig =
\sum_{i=0}^d w_i \, \alphag_i$ in $\s,$ and let $s=\sum_{i=0}^k w_i.$
We have
\begin{equation*}
\sum_{i=k+1}^d w_i=\sum_{i=0}^d w_i-\sum_{i=0}^k w_i=1-s;
\end{equation*}
and 
\begin{equation*}
\xig-F(\xig)=\sum_{i=0}^d w_i \, \alphag_i - \Frac{1}{s} \sum_{i=0}^k w_i \, \alphag_i
=\left(1-\frac{1}{s}\right)\left(\sum_{i=0}^k w_i\, \alphag_i \right) +
\sum_{i=k+1}^d w_i\, \alphag_i.
\end{equation*}
By the triangle inequality, we obtain
\begin{equation}\label{eq:norm-ineq}
|\xig-F(\xig)|\leq \max_{o \leq i \leq d} |\alphag_i| 
\left(\left| 1-\frac{1}{s} \right|\left(\sum_{i=0}^k w_i\right)  + 
\sum_{i=k+1}^d w_i\right)
=2 (1-s)\max_{0 \leq i \leq d} |\alphag_i|.
\end{equation}
Let $\l=\pi(\xig)=\sum_{i=k+1}^d w_i \, \pi(\alphag_i).$
Since $\pi(\alphag_i) \geq \l_0$ for all $i \geq k+1,$ it follows that
\begin{equation}
\l=\sum_{i=k+1}^d w_i \, \pi(\alphag_i) \geq \l_0 \left(\sum_{i=k+1}^d w_i
\right)= \l_0 (1-s).
\end{equation}
It follows that $1-s \leq \Frac{\l}{\l_0}.$ Combining this 
with~\eqref{eq:norm-ineq}, we obtain
\begin{equation*}
|\xig-F(\xig)|\leq 2 \,\frac{\l}{\l_0} \,\max_{0 \leq i \leq d} |\alphag_i| 
\leq 2 \,\frac{\l}{\l_0} \,\max\{ |\alphag|, \, \alphag \hbox{ vertex of }K\}.
\end{equation*}
Thus, $|\xig-F(\xig)|$ is bounded by a quantity independent of $\xig$
that goes to zero when $\l$ goes to zero, and since $|\x-f^\l(\x)|\leq
|\xig-F(\xig)|,$ the result follows.
\end{proof}

\begin{proposition}\label{prop:homeo}
For all $0<\l'<\l<\l_0,$ there exists a homeomorphism $h: X_\l
\to X_{\l'}$ such that $f^\l=f^{\l'}\circ h.$
\end{proposition}

\begin{proof}
Let $\xig \in X_\l,$ and $\D(\xig)$ be as in Proposition~\ref{prop:F}. 
Since $\pi$ is affine on $\D(\xig),$ if $\taug=t \xig+
(1-t)F(\xig)$ is a point on $\D(\xig),$ we have 
\begin{equation*}
\pi(\taug)=t\pi(\xig)+(1-t)\pi(F(\xig))=t \l.
\end{equation*}
Thus, $\taug\in X_{\l'}$ if and only if
$t=\l'/\l,$ and so the map $h$ defined by
\begin{equation}\label{eq:h}
h(\xig)=\frac{\l'}{\l}\,\xig + \left( 1- \frac{\l'}{\l}\right) F(\xig);
\end{equation}
maps $X_\l$ to $X_{\l'}.$

\medskip  

Suppose that there exists $\xig$ and $\xig'$ in $X_\l$ such that
$h(\xig)=h(\xig')=\taug.$ Then $\taug\in \D(\xig)$ and $\taug \in
\D(\xig'),$ and by Proposition~\ref{prop:F} this means that
$F(\xig)=F(\taug)=F(\xig').$ Then~\eqref{eq:h} implies that
$\xig=\xig',$ so $h$ is injective.  The map $h$ is also surjective,
since for $\taug \in X_{\l'},$ it is easy to verify that the point
$\xig$ defined by
\begin{equation*}
\xig=\frac{\l}{\l'}\taug-\left(\frac{\l}{\l'}-1\right)F(\taug);
\end{equation*}
is a point in $X_\l$ such that $h(\xig)=\taug.$

\medskip

The continuity of $h$ follows from the continuity of $F.$ Since $h(\xig)
\in \D(\xig)$ by construction, Proposition~\ref{prop:F} implies that
$F(h(\xig))=F(\xig),$ so $f^\l=f^{\l'} \circ h.$ 
\end{proof}

\section{Approximation of the fibered products}
\label{sec:approx}

We will now turn our attention to the fibered products associated to
the surjections $f^\l:X_\l \to X_0$ that were constructed in the
previous section. We will prove in Proposition~\ref{prop:approx} the
approximation result for those sets that yields the result of the main
theorem.

\medskip

Define for $p \in \N$ and $\l \in (0,\l_0),$ 
\begin{equation}\label{eq:W}
W_\l^p=\{(\x_0, \ldots, \x_p) \in (X_\l)^{p+1} \mid 
f^\l(\x_0)=\cdots=f^\l(\x_p)\}.
\end{equation}
{}From Theorem~\ref{th:ss}, 
we have for any $\l \in (0, \l_0),$
\begin{equation}\label{eq:ss}
b_k(X_0) \leq \sum_{p+q=k} b_q(W_\l^p).
\end{equation}
Thus, bounding the Betti numbers of $X_0$ can be reduced to estimating
the Betti numbers of the sets $W_\l^p$ for some $\l\in(0,\l_0).$ The
first step in that direction is the following.

\begin{proposition}\label{prop:Whomeo}
For all $0<\l'<\l<\l_0,$ the sets $W_\l^p$ and $W_{\l'}^p$ are
homeomorphic.
\end{proposition}

\begin{proof}
Fix $0<\l'<\l<\l_0,$ and let $h$ be the homeomorphism between $X_\l$
and $X_{\l'}$ described in Proposition~\ref{prop:homeo}.  Since
$f^{\l'}\circ h =f^\l,$ the map $h^p: (X_\l)^{p+1} \to (X_{\l'})^{p+1}$
defined by
\begin{equation}\label{eq:hp}
h^p(\x_0, \ldots, \x_p)=(h(\x_0), \ldots, h(\x_p));
\end{equation}
maps $W_\l^p$ homeomorphically onto $W_{\l'}^p.$
\end{proof}

Recall that for $p \in \N$ and $\x_0, \ldots,\x_p \in \R^n,$ $\rho_p$
is the polynomial
\begin{equation}
\rho_p(\x_0, \ldots,\x_p)=\sum_{0 \leq i < j \leq p} |\x_i-\x_j|^2.
\end{equation}
For $\l \in (0,\l_0),$  $\e>0$ and $\d>0,$ we define the following sets.
\begin{align*}
W_\l^p(\e)&=\{(\x_0, \ldots, \x_p) \in (X_\l)^{p+1} \mid 
\rho_p(f^\l(\x_0), \ldots, f^\l(\x_p)) \leq \e\};\\
D_\l^p(\d)&=\{(\x_0, \ldots, \x_p) \in (X_\l)^{p+1} \mid 
\rho_p(\x_0, \ldots, \x_p) \leq \d\}.
\end{align*}

\begin{proposition}\label{prop:Wehomeo}
Let $p \in \N$ be fixed. There exists $\e_0>0,$ such that for all $\l
\in (0,\l_0)$ and all $0<\e'<\e<\e_0,$ the inclusion $W_\l^p(\e') \inc
W_\l^p(\e)$ is a homotopy equivalence. In particular, this implies
that  
\begin{equation}\label{eq:bq}
b_q(W_\l^p(\e))=b_q(W_\l^p);
\end{equation}
for all integer $q,$  all $\l \in (0,\l_0)$ and all $\e\in(0,\e_0).$
\end{proposition}

\begin{proof}
First, notice that it is enough to prove the result for a fixed $\l
\in (0,\l_0),$ since if $0<\l'<\l<\l_0$ are fixed, the map $h^p$
introduced in~\eqref{eq:hp} induces a homeomorphism between
$W_{\l}^p(\e)$ and $W_{\l'}^p(\e)$ for any $\e>0.$

\medskip

Fix $\l \in (0,\l_0).$ 
By the generic triviality theorem (see \cite[Chapter~9,
Theorem~1.2]{vdd98} or \cite[Theorem~5.22]{coste}), there exists
$\e_0>0$ such that the projection
\begin{equation*}
\left\{ (\x_0, \ldots, \x_p, \e) \mid \e\in(0,\e_0) \hbox{ and }
(\x_0, \ldots, \x_p) \in W_\l(\e)\right\} \mapsto \e;
\end{equation*}
is a trivial fibration.
It follows that for all $0<\e'<\e<\e_0,$ the inclusion $W_\l^p(\e')
\inc W_\l^p(\e)$ is a homotopy equivalence, which proves the first
part of the proposition.

\medskip

The inclusions $W_\l^p(\e') \inc W_\l^p(\e)$
for $\e'<\e$ make the family $W_\l^p(\e), \e>0$ into a directed
system, and we have:
\begin{equation*}
\lim_{\longleftarrow} W_\l^p(\e) \iso \bigcap_{\e>0} W_\l^p(\e)=W_\l^p.
\end{equation*}
The induced maps in \v{C}ech homology make the family
$\ck{H}_*(W_\l^p(\e)), \e>0$ into a directed system too, and by the
continuity property of the \v{C}ech homology \cite[Chapter~10]{ES},
this implies that the \v{C}ech homology groups of $W^p_\l$ are
isomorphic to the inverse limit of the groups $\ck{H}_*(W_\l^p(\e)).$
Note however that the sets $W_\l^p$ and $W_\l^p(\e)$ being compact
definable sets, they are homeomorphic to finite simplicial
complexes, and thus, their singular and \v{C}ech homologies
coincide. Hence, we have for all $q,$
\begin{equation}\label{eq:limp}
H_q(W_\l^p) \iso \lim_{\longleftarrow} H_q(W_\l^p(\e)).
\end{equation}
Since the inclusion $W_\l^p(\e') \inc W_\l^p(\e)$ is a homotopy
equivalence for all $0<\e'<\e<\e_0,$ it induces an isomorphism in
homology for the directed system $H_q(W_\l^p(\e)).$ Thus, the ranks in
that system are constant and equal to the rank of the limit,
yielding~\eqref{eq:bq}.
\end{proof}

\begin{example}\label{ex:cex}
Consider the family $X\sub \R^2\times\R_+$ such that the fiber is
given for all $\l>0$ by
\begin{equation*}
X_\l=\{(x,0) \mid 0\leq x \leq 1\} \cup \{(0,\l)\} \cup
\{(\l,\l)\}. 
\end{equation*}
If we construct $f^\l$ from a triangulation as in
Definition~\ref{df:fl}, we must have $f^\l(0,\l)=f^\l(\l,\l)=0,$
independently of the choice of the triangulation.  But there are other
families of continuous surjections from $X_\l$ into $X_0$ that are
close to identity. A natural choice would be for instance to take the
maps $g^\l$ defined by $g^\l(x,y)=(x,0),$ for which we still have
\begin{equation*}
\lim_{\l\to 0} \max_{(x,y)\in X_\l} |g^\l(x,y)-(x,y)|=0.
\end{equation*}
However, it's easy to check that the topological type of the
corresponding sets $W^2_\l(\e)$ changes exactly when $\e=\l,$ since
the two connected components of $(X_\l)^2$ formed by the isolated
points $\{(0,\l,\l,\l)\}$ and $\{(\l,\l,0,\l)\}$ belong to
$W^2_\l(\e)$ exactly when $\e\geq \l.$
\end{example}

Example~\ref{ex:cex} shows that finding a family of surjections close
to identity is not enough to guarantee the existence of an $\e_0$ {\em
independent of} $\l$ for which Proposition~\ref{prop:Wehomeo}
holds. This explains why a careful construction of $f^\l$ was
necessary in section~\ref{sec:triang}.  The fact that we can find
$\e_0$ independent of $\l$ will play a key part to approximate the
sets $W^p_\l$ (see Proposition~\ref{prop:approx}).

\bigskip

\begin{proposition}\label{prop:deltai}
Let $p\in \N$ be fixed. 
There exists $\l_1$ such that $0<\l_1\leq \l_0$ and definable
functions $\d_0(\l)$ and $\d_1(\l)$ defined for $\l\in(0,\l_1)$ such
that $\lim_{\l \to 0} \d_0(\l)=0,$ $\lim_{\l \to 0} \d_1(\l)\neq0,$
and such that for all $0<\d_0(\l)< \d'<\d<\d_1(\l),$ the inclusion
$D_\l^p(\d') \inc D_\l^p(\d)$ is a homotopy equivalence.
\end{proposition}

\begin{proof}
Let $\l\in(0,\l_0)$ be fixed.
As in the proof of Proposition~\ref{prop:Wehomeo}, we can apply the
generic triviality theorem to the projection of the family
$\{D^p_\l(\d) \mid \d>0\}$ on the $\d$-axis. This yields 
real numbers $d_0(\l)=0< d_1(\l) < \cdots< d_m(\l) <
d_{m+1}(\l)=\infty$ such that for all $0\leq i\leq m$
the projection is a trivial fibration above the interval
$(d_i(\l),d_{i+1}(\l)).$ 

\medskip

We will show now that we can choose the numbers $d_i(\l)$ to depend
definably on $\l$ for $\l$ in an interval of the from $(0,\l_1).$ 
Consider the definable set
\begin{equation*}
\Dcal=\{(\x_0, \ldots, \x_p, \l, \d) \mid 
(\x_0, \ldots, \x_p) \in D^p_\l(\d), \l>0, \d>0\},
\end{equation*}
and its projection on the $(\l,\d)$-plane. By the generic triviality
theorem, the region $\{\l>0,\d>0\}$ can be partitioned in finitely
many definable subsets such that the projection of $\Dcal$ is a
trivial fibration over each of those subsets. Without loss of
generality, we can assume that those subsets are cylindrical cells for
the $\l$-projection (this is the o-minimal equivalent of cylindrical
algebraic decomposition, see \cite[Chapter~3]{vdd98} or
\cite[Definition~2.4]{coste}). In particular, there exists $\l_1>0$
and continuous definable functions $d_0(\l)=0< d_1(\l) < \cdots<
d_m(\l) < d_{m+1}(\l)=\infty$ defined for $\l\in(0,\l_1)$ such that
the cells
\begin{equation*}
C_i=\{(\l,\d) \mid \l\in(0,\l_1), d_i(\l) < \d < d_{i+1}(\l)\},
\end{equation*}
are part of this partition for $0 \leq i \leq m.$ Without loss of
generality, we can assume that $\l_1 \leq \l_0.$

\medskip

The functions $d_i(\l)$ being definable, each has a well-defined
(although possibly infinite) limit when $\l$ goes to zero.  Let $j$ be
the largest index such that $\lim_{\l \to 0} d_j(\l)=0,$ and 
let $\d_0(\l)=d_j(\l)$ and $\d_1(\l)=d_{j+1}(\l).$ The functions
$\d_0$ and $\d_1$ are definable and verify
\begin{equation*}
\lim_{\l \to 0} \d_0(\l)=0 \hbox{ and } \lim_{\l \to 0} \d_1(\l)>0.
\end{equation*}
Since the projection of $\Dcal$ is a trivial fibration over the cell
$C_j,$ then for any fixed $\l \in (0,\l_1)$ and any $0\leq\d_0(\l)<
\d'<\d<\d_1(\l),$ the inclusion $D_\l^p(\d') \inc D_\l^p(\d)$ is
certainly a homotopy equivalence.
\end{proof}

\bigskip

Let $R>0$ be such that $X_\l \sub \{|\x| \leq R\}$ for all $\l \in
(0,1).$ We define for $p \in \N,$ 
\begin{equation}\label{eq:etap}
\eta_p(\l)=p(p+1)\,\left(4R\,\max_{\x\in X_\l} |\x-f^\l(\x)|
+2\,(\max_{\x\in X_\l} |\x-f^\l(\x)|)^2\right).
\end{equation}
By Proposition~\ref{prop:normF}, we have 
\begin{equation*}
\lim_{\l \to 0} \eta_p(\l)=0.
\end{equation*}

\begin{lemma}\label{lem:inclusions}
For all $\l \in (0,\l_0),$ $\d>0$ and $\e>0,$ the following inclusions
hold.
\begin{equation*}
D_\l^p(\d) \sub W_\l^p(\d+\eta_p(\l)), \hbox{ and } 
W_\l^p(\e) \sub D_\l^p(\e+\eta_p(\l)).
\end{equation*}
\end{lemma}

\begin{proof}
Let $m(\l)=\max_{\x\in X_\l} |\x-f^\l(\x)|.$
For any $\x_i, \x_j$ in $X_\l,$ the triangle inequality gives
\begin{align*}
|f^\l(\x_i) - f^\l(\x_j)|^2 &\leq [|f^\l(\x_i)-\x_i| + |\x_i - \x_j| +
|\x_j - f^\l(\x_j)|]^2\\
&\leq [|\x_i-\x_j|+2m(\l)]^2\\
&\leq |\x_i-\x_j|^2 + 8R \, m(\l)+
4 m(\l)^2.
\end{align*}
Summing this inequality for all $0 \leq i < j \leq p,$ we obtain that
for any $\x_0, \ldots, \x_p$ in $X_\l,$
\begin{equation*}
\rho_p(f^\l(\x_0), \ldots, f^\l(\x_p)) \leq \rho_p(\x_0, \ldots,
\x_p)+\eta_p(\l). 
\end{equation*}
The first inclusion follows easily from this inequality. The
second inclusion follows from a similar reasoning.
\hbox{ }
\end{proof}

\begin{proposition}\label{prop:approx}
For any $p\in \N,$ there exists $\l \in (0,\l_0),$ $\e\in (0,\e_0)$
and $\d>0$  such that
\begin{equation}
H_*(W_\l^p(\e)) \cong H_*(D_\l^p(\d)).
\end{equation}
\end{proposition}

\begin{proof}
Let $\d_0(\l)$ and $\d_1(\l)$ be the functions defined in
Proposition~\ref{prop:deltai}. Since the limit when $\l$ goes to zero
of $\d_1(\l)-\d_0(\l)$ is not zero, whereas the limit of $\eta_p(\l)$
is zero, we can choose $\l>0$ such that
$\d_1(\l)-\d_0(\l)>2\eta_p(\l).$ Then, we can choose $\d'>0$ such that
$\d_0(\l)<\d'<\d'+2\eta_p(\l)< \d_1(\l).$ Taking a smaller $\l$
if necessary, we can also assume that $\d'+3\eta_p(\l)<\e_0.$

\medskip

Let $\e=\d'+\eta_p(\l),$ $\d=\d'+2\eta_p(\l)$ and
$\e'=\d'+3\eta_p(\l).$ From Lemma~\ref{lem:inclusions}, we have the
following sequence of inclusions;
\begin{equation*}
D_\l^p(\d')\stackrel{i}{\inc} W_\l^p(\e)
\stackrel{j}{\inc} D_\l^p(\d) \stackrel{k}{\inc} W_\l^p(\e').
\end{equation*}

Since both $\e$ and $\e'$ are less than $\e_0,$
Proposition~\ref{prop:Wehomeo} ensures that the inclusion map $k \circ
j$ is a homotopy equivalence.  Similarly, since both $\d$ and $\d'$
are in the interval $(\d_0(\l),\d_1(\l)),$ it follows from
Proposition~\ref{prop:deltai} that the inclusion $j \circ i$ is a
homotopy equivalence too. In particular, both inclusions give rise to
isomorphisms on the homology level.  The resulting diagram in homology
is the following;
\begin{equation*}
\xymatrix{
H_*(D_\l^p(\d'))\ar[rr]^{(j\circ i)_*}_{\cong} 
\ar[dr]^{i_*}
& & H_*(D_\l^p(\d))\ar[dr]^{k_*} & \\
& H_*(W_\l^p(\e))\ar[rr]^{(k\circ j)_*}_{\cong} 
\ar[ur]^{j_*}
& & H_*(W_\l^p(\e'))
}
\end{equation*}
Since $(j \circ i)_*=j_*\circ i_*,$ the surjectivity of $(j \circ
i)_*$ implies that $j_*$ must be surjective, and similarly, the fact
that $(k \circ j)_*=k_*\circ j_*$ is injective implies that $j_*$ is
injective. Hence, $j_*$ is an isomorphism between $H_*(W_\l(\e))$ and
$H_*(D_\l(\d)),$ as required.
\end{proof}

\bigskip
\noindent
\textbf{Proof of Theorem~\ref{th:main}.} The proof of the main theorem
follows now easily from the results in this section. If $L$ is the
Hausdorff limit of a sequence of fibers $A_{a_i},$ we can construct a
family $X$ as in Proposition~\ref{prop:red} such that $X_0=L.$ Then,
for all $\l \in (0, \l_0)$, we have 
\begin{equation*}
b_k(L) \leq \sum_{p+q=k} b_q(W^p_\l).
\end{equation*}
{}From Proposition~\ref{prop:approx}, for $\l$ small enough, we can find
$\e\in (0,\e_0)$ and $\d>0$ such that
$b_q(W^p_\l)=b_q(W^p_\l(\e))=b_q(D^p_\l(\d))$ for every integer $0
\leq p \leq k.$ Thus, if $a\in A'$ is such that 
$X_\l=A_a,$  inequality~\eqref{eq:bkL} in
Theorem~\ref{th:main} holds for that $a.$ \hfill $\Box$

\begin{rem}
Note that the upper-bound in Theorem~\ref{th:main} depends only on the
sets $D^p_a(\d),$ which are defined from the original family $A$ but
independent of the auxiliary family $X.$ Thus, we'll be able to
derive quantitative estimates from our main result without ever having
to worry about the complexity of the fibers $X_\l.$
\end{rem}

This is an important remark, because the complexity of the description
of $X_\l$ may be much worse than the complexity of a fiber $A_a,$ 
even for a choice of $a$ such that the two sets are equal.

\section{Effective estimates on the Betti numbers}
\label{sec:estimates}

This section is devoted to the quantitative estimates than can be
derived from Theorem~\ref{th:main}, both in the algebraic and the
Pfaffian case.

\subsection{Pfaffian functions}

We'll start by recalling the basic results about Pfaffian functions.
Let $\U \sub \R^n$ be an open domain. 

\begin{definition}
Let $\x=(x_1, \ldots, x_n)$ and let $(f_1(\x), \ldots, f_\ell(\x))$ be
a sequence of analytic functions in $\U.$ This sequence is called a
{\em Pfaffian chain} if the functions $f_i$ are solution on $\U$ of a
triangular differential system of the form;
\begin{equation}\label{eq:chain}
df_i(\x)=\sum_{j=1}^n P_{i,j}(\x, f_1(\x), \ldots, f_i(\x))dx_j;
\end{equation}
where the functions $P_{i,j}$ are polynomials in $\x, f_1, \ldots,
f_i.$

If $(f_1, \ldots, f_\ell)$ is a fixed Pfaffian chain on a domain $\U,$
the function $q$ is a {\em Pfaffian function} in the chain $(f_1,
\ldots, f_\ell)$ if there exists a polynomial $Q$ such that for all $\x
\in \U,$ 
\begin{equation}\label{eq:Q}
q(\x)=Q(\x,f_1(\x), \ldots, f_\ell(\x)).
\end{equation}
\end{definition}

\medskip

Pfaffian functions come naturally with a notion of complexity, or {\em
format.}  If $(f_1, \ldots, f_\ell)$ is a Pfaffian chain, we call
$\ell$ its {\em length}, and we let its {\em degree} $\alpha$ be the
maximum of the degrees of the polynomials $P_{i,j}$ appearing
in~\eqref{eq:chain}.  If $q$ is as in~\eqref{eq:Q}, the degree $\beta$
of the polynomial $Q$ is called the {\em degree} of $q$ in the chain
$(f_1, \ldots, f_\ell).$ 

\begin{definition}\label{df:qformat}
For $q$ as above, the tuple $(n,\ell,\alpha,\beta)$ is called the {\em
format} of $q.$
\end{definition}

\begin{example}\label{ex:few}
Let $\x=(x_1, \ldots, x_n) \in \R^n$ and $\m_1, \ldots, \m_\ell$ be
fixed vectors in $\R^n.$ Define for all $1 \leq i \leq \ell,$
$f_i(\x)=e^{\langle \m_i, \x\rangle},$ (where
$\langle\cdot,\cdot\rangle$ denotes the Euclidean scalar product).
Then, $(f_1, \ldots, f_\ell)$ is a Pfaffian chain of length $\ell$ and
degree $\alpha=1$ on $\R^n.$
\end{example}

The class of Pfaffian functions is a very large class that contains,
among other things, all real elementary functions (for a suitable
choice of the domain of definition), Liouville functions, and Abelian
integrals.  We refer the reader to the book~\cite{k91} or the
papers~\cite{gv:strat,gv:survey} for more details and examples.  The
main result about Pfaffian functions is the following.

\begin{theorem}[Khovanskii~\cite{k91}]\label{th:khov}
Let $(f_1, \ldots, f_\ell)$ be a Pfaffian chain of length $\ell$ and
degree $\alpha$ defined on $\R^n.$ Let $(q_1, \ldots, q_n)$ be Pfaffian
functions in that chain, and let $(n,\ell,\alpha,\beta_i)$ be the
format of $q_i$ for $1 \leq i \leq n.$
Then the number of solutions of the system
\begin{equation}\label{eq:thkh}
q_1(\x)=\cdots=q_n(\x)=0,
\end{equation}
that are isolated in $\C^n$ is bounded from above by
\begin{equation}\label{eq:kh}
2^{\ell(\ell-1)/2} \, \beta_1\cdots\beta_n\; (\beta_1+\cdots+\beta_n
-n+\min(n,\ell)\alpha+1)^\ell.
\end{equation}
\end{theorem}

\begin{rem}
In particular, using Example~\ref{ex:few} through a logarithmic change
of variables, one can show that if $(q_1, \ldots, q_n)$ are sparse
real polynomials, the number of isolated roots of the system
$q_1(\x)=\cdots=q_n(\x)=0$ in the quadrant $(\R_+)^n$ can be bounded
independently of the degrees of the polynomials $q_i.$ If no more than
$\ell$ monomials appear with a non-zero coefficient in at least one of
the polynomials, Theorem~\ref{th:khov} gives that the number of roots
of the system is bounded by $2^{\ell(\ell-1)/2} \, (n+1)^\ell.$
(See~\cite{k91}. Note that this bound is known to be pessimistic, at
least in some cases, see~\cite{lirojaswang} for instance.)
\end{rem}

Theorem~\ref{th:khov} can be easily generalized to the case where the
Pfaffian chain $(f_1, \ldots, f_\ell)$ is defined on a domain $\U$
different from $\R^n.$ In that case, the bound~\eqref{eq:kh} becomes
\begin{equation}\label{eq:khgal}
2^{\ell(\ell-1)/2} \, \beta^n\, O(n(\alpha+\beta))^\ell.
\end{equation}
The constant coming from the $O(\cdots)$ notation depends
only on the geometry of the open domain $\U.$ In many cases (for
instance if the functions $q_i$ are real elementary functions), the
domain $\U$ is given by 
\begin{equation*}
\U=\{\x \mid g_1(\x)>0, \ldots, g_r(\x)>0\};
\end{equation*}
where the functions $g_i$ are Pfaffian. In that case, a suitable
constant can be determined explicitly (see for instance Chapter~1 in
\cite{z:these} for a discussion of the determination of such constants).

\subsection{Quantifier-free formulas}
In this section, we will let $\P=\{p_1, \ldots, p_s\}$ be either a set
of polynomials, or a set of Pfaffian functions defined in a common
Pfaffian chain $(f_1, \ldots, f_\ell)$ on a definable domain $\U$. By a
quantifier-free formula on $\P,$ we mean a Boolean combination of sign
conditions on the functions in $\P,$ as defined below.

\begin{definition}\label{def:qf}
A formula $\Phi$ is called a {\em quantifier-free formula on $\P$} if it is
derived from atoms of the form $p_i \star 0,$ -- where $1 \leq i \leq
s$ and $\star \in \{=,\leq,\geq\},$ -- using conjunctions,
disjunctions and negations. Moreover, we will say that the formula
$\Phi$ is {\em $\P$-closed } if it was derived without using negations.
\end{definition}

We endow quantifier-free formulas with the following
format.\footnote{%
Note that this may not be the most standard notion of format for
quantifier-free formulas, but it is well-adapted to the Betti number
bounds we will be discussing next.}

\begin{definition}\label{format}
Let $\Phi$ be a quantifier-free formula on $\P,$ where $\P$ is a
collection of $s$ polynomials in $n$ variables of degree bounded by
$d.$ Then $(n,d,s)$ is called the {\em format} of $\Phi.$

Similarly, if $(f_1, \ldots, f_\ell)$ is a fixed Pfaffian chain and
$\P$ is a collection of $s$ functions in that chain, where the format
of those functions is at most $(n,\ell,\alpha,\beta),$
then the {\em format} of any quantifier-free formula on $\P$ is  
$(n,\ell, \alpha, \beta, s).$ 
\end{definition}

If $\Phi$ is a quantifier-free formula on $\P$, where $\P$ is a
collection of polynomials, the associated {\em semialgebraic} set is
$X=\{\x \in \R^n \mid \Phi(\x)\}.$ If $\P$ is a collection of Pfaffian
functions in a Pfaffian chain defined on a domain $\U$ and $\Phi$ is a
quantifier-free formula on $\P$, the associated {\em semi-Pfaffian
set} is the set $X=\{\x \in \U \mid \Phi(\x)\}$.

\medskip

Both polynomials and Pfaffian functions generate o-minimal
structures. For polynomials, this is a consequence of the famous
Tarski-Seidenberg theorem (see~\cite{bcr} or~\cite{vdd98}).  The
o-minimality of the structure $\S_\pfaff$ generated by Pfaffian
functions was first proved by Wilkie~\cite{wilkie99}, and generalized
in~\cite{KaMcI,lr98,s99}. Gabrielov introduced the notion of {\em
relative closure} in~\cite{g00} to offer a description of $\S_\pfaff$
more adapted to solving quantitative questions in that structure,
and  Corollary~\ref{cor:rcbnd} is an example of how this construction
can yield answers in the case of topological complexity (see
also~\cite{gz} for estimates on the number of connected components).
The relation between relative closures and Hausdorff limits 
is discussed in more details in the introduction. The
reader can also refer to~\cite{g00,gv:survey,gz} for a precise
description of the construction of $\S_\pfaff$ via relative closures.

\subsection{Betti number bounds for quantifier-free formulas}

Elementary Morse theory used in conjunction with the B\'ezout
inequality lets us bound the Betti numbers of real algebraic varieties
in terms of the degrees of the defining polynomials. This was first
noticed by Oleinik and Petrovsky, and then independently by Thom and
Milnor (the reader can refer to \cite[Chapitre~11]{bcr} for a proof
and the corresponding references). Khovanskii noted that the method
could be generalized to zero-sets of Pfaffian functions~\cite{k91,
z99}.

\bigskip

The notion of $\P$-closed formulas was introduced by Basu in the
algebraic setting. Using basic algebraic topology techniques and 
the bound of Oleinik, Petrovsky, Milnor and Thom, he was able to prove
the following estimate.%
\footnote{Basu's estimate is actually slightly sharper than the one
given here, but the simpler form~\eqref{eq:aPclosed} is enough 
for our purposes.}

\begin{theorem}[\cite{basu}]\label{th:aPclosed}
Let $\Phi$ be a $\P$-closed polynomial formula of format $(n,d,s)$ and
let $X$ be the corresponding semialgebraic set.  The sum of Betti
numbers of $X$ admits the upper-bound;
\begin{equation}\label{eq:aPclosed}
b(X) \leq  O(sd)^n.
\end{equation}
\end{theorem}

Basu's technique goes through in the Pfaffian setting, and one obtains
the following result.

\begin{theorem}[\cite{z99}]\label{th:Pclosed}
Let $(f_1, \ldots, f_\ell)$ be a Pfaffian chain  on a definable domain
$\U.$ Let $\P$ be a collection of Pfaffian functions in that chain,
and let $\Phi$ be a $\P$-closed formula of format $(n,\ell, \alpha,
\beta, s).$ If $X$ is the corresponding semi-Pfaffian set $X=\{\x
\in \U \mid \Phi(\x)\},$ the sum of the Betti numbers of $X$ admits
the following upper-bound;
\begin{equation}\label{eq:Pclosed}
b(X) \leq 2^{\ell(\ell-1)/2} s^n O(n (\alpha+ \beta))^{n+\ell};
\end{equation}
where the constant depends only on the definable domain $\U.$
\end{theorem}

Note that the above bound requires that the constant depends on $\U$
because Khovanskii's theorem applied on a general domain $\U$ gives a
bound of the form~\eqref{eq:khgal}, that involves a constant depending
on the geometry of $\U.$

\begin{rem}
For any estimate on the Betti numbers of semi-Pfaffian sets, we must
assume that the domain $\U$ under consideration is definable.  Indeed,
one can easily construct non-definable domains $\U$ for which there
exists semi-Pfaffian sets $X$ such that $b(X)$ is infinite.
\end{rem}

In the case where $X$ is compact, but is {\em not} given by a
$\P$-closed formula (for instance, if $X$ is a semialgebraic set that
has been obtained by a quantifier elimination procedure), bounds on
the Betti numbers can be established using the fact that the singular
homology of $X$ and its Borel-Moore homology coincide when $X$ is
compact. Using the sub-additivity of the Borel-Moore homology (see for
instance \cite[Proposition~1.8]{momopa} \cite{yao} \cite{burgisser})
along with Theorem~\ref{th:aPclosed} and the estimates on the number
of sign cells in~\cite{bpr:cells}, one can show that the rank of the
Borel-Moore groups of a compact semialgebraic set defined by a
formula of format $(n,d,s)$ is bounded by $O(sd)^{2n}.$ The same
method yields an estimate for semi-Pfaffian sets
too~\cite[Theorem~2.23]{z:these}. However, these results have been
superseded by Theorem~\ref{th:aqf} below.

\bigskip

Recently, Gabrielov and Vorobjov used Theorem~\ref{th:aPclosed} to
establish a Betti number bound valid for {\em any} semialgebraic set
$X,$ with no assumptions either on the topology of $X$ or on the shape
of its defining formula. The result they obtained is the following.

\begin{theorem}[\cite{gv:qf}]\label{th:aqf}
Let $\Phi$ be a polynomial quantifier-free formula of format
$(n,d,s)$ and let $X$ be the corresponding semialgebraic set.
The sum of Betti numbers of $X$ admits the following upper-bound.
\begin{equation}\label{eq:aqf}
b(X) \leq O(s^2d)^n.
\end{equation}
\end{theorem}

Again, the result can be generalized without problem to the Pfaffian
case, to obtain the following estimate.

\begin{theorem}\label{th:qf}
Let $X$ be {\em any} semi-Pfaffian set defined by a quantifier-free
formula of format $(n,\ell, \alpha, \beta, s).$ 
The sum of the Betti numbers of $X$ admits a bound of the
form 
\begin{equation}\label{eq:qf}
b(X) \leq 2^{\ell(\ell-1)/2} s^{2n} O(n (\alpha+ \beta))^{n+\ell};
\end{equation}
where the constant depends only on the definable domain $\U.$
\end{theorem}

\subsection{Proof of Corollary~\ref{cor:algebraic} %
and Corollary~\ref{cor:rcbnd}}

Let $A \sub \R^n \times \R^r$ be a definable family of compact,
uniformly bounded sets defined by a quantifier-free formula $\Phi(\x,a)$
whose atoms are sign condition on either polynomials or Pfaffian
functions, and let $L$ be a Hausdorff limit of fibers in $A.$ 
According to Theorem~\ref{th:main}, we
can bound $b_k(L)$ by estimating $b(D_a^p(\d))$ for all $0 \leq p
\leq k$ and suitable values $a=a^*$ and $\d=\d^*$ (the precise values of
$a^*$ and $\d^*$ do not affect the estimate).

\medskip

Suppose first that the atoms of $\Phi(\x,a)$ are polynomials. Then
for any $p,\d,a,$ the set $D_a^p(\d)$ is given by the
quantifier-free formula $\Psi_p(\x_0,\ldots,\x_p,a)$ such that
\begin{equation}\label{eq:Dqf}
\Psi_p(\x_0,\ldots,\x_p,a)=
\Phi(\x_0,a) \wedge \cdots \wedge \Phi(\x_p,a) \wedge \rho_p(\x_0,
\ldots, \x_p)\leq \d.
\end{equation}
The set of functions $\Q_p$ appearing in the atoms of $\Psi_p$ are of the
form $f(\x_i,a),$ for any $0 \leq i \leq p$ and any $f\in \P,$ plus
of course $\rho_p(\x_0, \ldots, \x_p)-\d.$ Thus, specializing $\Psi_p$
at $a=a^*$ gives a formula of format $(n(p+1),d,s(p+1)+1),$ where
$s=\P$ and
\begin{equation*}
d=\max(2, \{\deg_\x f(\x,a) \mid f\in \P\}).
\end{equation*}
Theorem~\ref{th:aqf} gives $b_q(D_a^p(\d)) \leq O(p^2s^2d)^{(p+1)n},$
and Corollary~\ref{cor:algebraic} follows.

\medskip

The same reasoning goes through in the Pfaffian case. If after
specialization, the format of $\Phi(\x,a^*)$ is $(n,\ell, \alpha,
\beta, s),$ the format of $\Psi_p(\x_0,\ldots,\x_p,a^*)$ is 
$(n(p+1),\ell(p+1),\alpha, \beta, s(p+1)+1).$
\footnote{The length of the Pfaffian chain is multiplied by $p+1$
since we need a copy of each function in the chain for every one of
the blocks of variables $\x_0, \ldots,\x_p.$}
From Theorem~\ref{th:qf}, we obtain that
\begin{equation*}
b(D_a^p(\d)) \leq 2^{\ell(p+1)[\ell(p+1)-1]/2} \, s^{2n(p+1)} \,
O(np(\alpha+\beta))^{(p+1)(n+\ell)},
\end{equation*}
and thus~\eqref{eq:H-betti} holds for $L.$
This is true for any number of parameters $r$ in the family $A,$ and
Corollary~\ref{cor:rcbnd} is simply the special case where $r=1.$
\hfill $\Box$

\begin{rem}
If the fiber $A_a$ is defined by a $\P$-closed formula, then the
formula $\Psi_p$ defining $D_\l^p(\d)$ is a $\Q_p$-closed formula, and
the dependence on $s$ of the bound on $b_k(L)$ can be improved from
$s^{2n(k+1)}$ to $s^{n(k+1)}$ by using the tighter estimates available
in that case (Theorem~\ref{th:aPclosed} in the algebraic case and
Theorem~\ref{th:Pclosed} in the Pfaffian case).
\end{rem}

\begin{acknowledgements}
I am grateful to my advisors Andrei Gabrielov and Marie-Fran\c{c}oise
Roy for their help and support. I also wish to thank Saugata Basu,
Michel Coste, Jean-Marie Lion and Nicolai Vorobjov for helpful
discussions, and the anonymous referees for their careful reading and
their many pertinent suggestions.
\end{acknowledgements}

\end{document}